\documentclass[12pt, twoside,a4paper]{article}
\pdfoutput=1
\usepackage[cp1251]{inputenc}
\usepackage{ifthen}
\usepackage{amsmath}
\usepackage{amsfonts}
\usepackage{latexsym}
\usepackage{amsthm}
\usepackage{amssymb}
\newcommand{\CopyName}{ V.\ M.\ Zhuravlov}
\newcommand{\NAME}{ V.\ M.\ Zhuravlov}
\newcommand{\Year}{2024}
\newcommand{\rightheadtext}{Nonclassical logics and multivariate truth values}
     \newcounter{chapter}
     \newcounter{artpage}[chapter]
     \newcommand{\vs}{\vspace{.1in}}
     \newcommand{\vsk}{\vspace{.2in}}
\makeatletter
     \renewcommand{\@evenhead}{\footnotesize \ifthenelse{\value{artpage}=0}
     {\hfil}{\thepage\hfil \textsc {\leftmark} \hfil } }
     \renewcommand{\@oddhead}{\footnotesize\ifthenelse{\value{artpage}=0}
     {\hfil}{\hfil \textsc \rightmark \hfil \thepage} }
     \newcommand{\logo}{\baselineskip2pc \hbox to\hsize{\hfil\copyright\,\footnotesize
     \CopyName, \Year}}
     \renewcommand{\@oddfoot}{\ifthenelse{\value{artpage}=0}{\logo
     \refstepcounter{artpage}} {\hfil\refstepcounter{artpage}}}
     \renewcommand{\@evenfoot}{\ifthenelse{\value{artpage}=0}{\logo
     \refstepcounter{artpage}} {\hfil\refstepcounter{artpage}}}
     \renewcommand{\section}{\@startsection{section}{1}{0pt}{3.5ex plus
     1ex minus .2ex}{2.3ex plus 2.ex}{\large\hfil\textsc}}
     \vspace{3pt}\vs

\vs
\newcommand{\tit}{Nonclassical logics and multivariate truth values}
\date{2023}
\usepackage{hyperref}
\usepackage{graphicx}
\graphicspath{{Images/}}
\begin{document}
\hfill
\vspace{0.3in}
\markboth{{\NAME}}{{\rightheadtext}}\begin{center} \textsc {\CopyName} \end{center}\begin{center} \renewcommand{\baselinestretch}{1.3}\bf {\tit} \end{center}
\vspace{20pt plus 0.5pt} {\abstract{\noindent
The article demonstrates that logic is not necessarily singleton and does not always have the standard interpretation of negation. Appropriate generalizations of logic are suggested. Positive logic and multivalued negation operations are examined. An attempt is made to abstract from the fixed and absolute values of truth and falsehood; the ambiguity of logical values is interpreted as their relativity,— projective logic\newline
\textit{Nonclassical logics and multivariate truth values, 2024, msc: 03H99;\vspace{3pt}}\newline
}
}\vsk
\section{Introduction}\par
Many types of non-classical logics are currently being studied. And most often the reason for studying them is their practical relevance. But what is practice? What is in demand today was not so yesterday and may not be so tomorrow. From this point of view, non-classical logics are a cultural phenomenon, despite all their mathematical rigor. However, this state of affairs can be observed not only in logic, but also in some other areas of mathematics, for example, in mathematical statistics or in economics; these are all scientific fields that do not yet have the empirical transparency and rigor of physics or chemistry. As for logic, it would be interesting to look more systematically at its non-classical generalizations in order to identify the main trends.\newline
\section{Cyclic many-valued logics with varying truth}\par
Let us turn our attention to ternary cyclic orders. Axioms of complete ternary order:\par
cyclicity: $[a,b,c]\Longrightarrow[b,c,a]$\par
antisymmetry: $[a,b,c]\Longrightarrow\neg[c,b,a]$\par
transitivity: $([a,b,c]\bigwedge[a,c,d])\Longrightarrow\neg[a,b,d]$\par
completeness: $((a\neq b)\bigwedge(b\neq c)\bigwedge(c\neq a))\Longrightarrow[a,b,c]\bigvee[c,b,a]$\par
This ternary relation is obtained by localization from the usual binary linear ordering. And this localization is consistent with respect to the elements: if we exclude the axiom of cyclicity, we obtain a single global linear order. And if we fix an arbitrary element a, then for each such element we obtain our own linear ordering on the entire set of values $(b\leq c)\equiv[a,b,c]$.
By non-singletonic logic we mean logic that is abstracted from the strict division of sets into elements.\par
Our idea here is that by choosing an arbitrary a as, say, a false value, we obtain our own system of many-valued logic; the result will be not just a many-valued logic, but a logic with flexible truth values (these can be different values for different elements of the model set — in fact, we will not even get scalar truth values, but vector ones).\newline
\section{Vector logic}\par
Every Boolean algebra is a Boolean ring with respect to the symmetric difference (as a ring sum) and the conjunction as multiplication. And, moreover, we can take an arbitrary element of this ring and create a 2-element Boolean algebra from it; then add another such Boolean algebra, creating a minimal algebra from a pair of such elements; continuing this process, if necessary, a transfinite number of times, we will represent the original algebra as the direct product of a certain set of such 2-element rings, obtaining a vector space (with a fixed basis) over a 2-element ring.\par
However, a 2-element ring consists of only zero and one and is also a field (each element has an inverse with respect to both addition and multiplication). Thus, every Boolean algebra is a vector space with a fixed basis over a 2-element field. By defining other bases in this space, we obtain a certain (not necessarily complete) system of subalgebras; and by embedding the original 2-element field into its extensions (up to real or even complex fields), we can reach the complete system of subalgebras of the original algebra. Then we get that the original algebra is a subalgebra of a vector space over a real (or complex) field. In general, this space has a transfinite basis. Moreover, in such a space one can study objects that are more "flexible" than just elements of the original Boolean algebra (linear combinations of different elements or sets).\par
Thus we arrive at the vector spaces studied by modern physics. What is important for us now is that quantum logic leads to lattices of vector subspaces that possess the property of modularity. Later we will encounter exactly such lattices when describing logic.\newline
\section{Non-singleton logic}\par
By non-singletonic logic we mean logic that is abstracted from the strict division of sets into elements.\par
The internal logic of a topos is defined by a subobject classifier together with some (logical) operations on it. A topos is a categorical generalization of set theory. A topos is a Cartesian closed category with a subobject classifier. A subobject classifier is a generalization of the logic of set theory.\par
It is important for us to generalize definition of subobject classifier so as to eliminate elements (i.e., singletons) from it. The subobject classifier is defined up to isomorphism. That is, there may be more than one element that can play the role of "truth".\begin{center}
\textbf{Generalized subobject classifier}\end{center}
~~~~Trying to find generalizations of logic that abstract away from elements, we should arrive at a true function, with which we can uniquely represent the algebra of subobjects of any object without using arrows from the terminal object. All this leads to the idea of some universal subobject. Thus we obtain a certain modification of the axiom of the existence of a subobject classifier:\par
\textbf{Axiom} \textit{S}: There is an object \textbf{\textit{D}}, together with the arrow $true:S\rightarrow D$, such that for every monoarrow $f:b\rightarrow a$ there is a unique pair of arrows $x_f:a\rightarrow D, y_f:b\rightarrow S$, for which there will be a Cartesian square:\begin{center}
\includegraphics{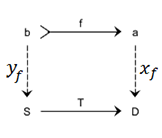}\end{center}
We will take this statement as an axiom of the subobject classifier. Then the natural definition would be:\par
\textbf{Definition 1}: If a category has inverse images, amalgams and satisfies the axiom \textbf{\textit{S}}, we will call it a logical category.\par
That is, unlike a topos, we are not talking here about completeness, or cocompleteness, exponentiation, and terminal or initial objects, — since we consider these axioms of a topos to be irrelevant to logic (they are relevant to set theory).\par
\textbf{Important notes.}\par
We can formulate the \textbf{\textit{S}}-axiom as follows: For each subobject there is a unique arrow $x_f:a\rightarrow D$, for which the corresponding Cartesian square. This is a weaker option. And I don't know which one is better.
Further, the existence of amalgams is not at all necessary for our research. However, it can be useful for defining some logical functions (for example, the supremum of a pair of statements; but one can also talk about logic without the “or” operator).
In addition, object \textbf{\textit{S}} may coincide with object \textbf{\textit{D}}, — $true:D\rightarrow D$. The meaning of this coincidence is that any arrow $D\rightarrow D$ replaces a certain fixed element \textbf{\textit{D}} with another that depends on it. This indicates not only some flexibility in fixing elements, but also the possibility of their localization relative to each other.
But, unless otherwise stated, we will adhere to the general \textbf{definition 1}.\par
Let us define truth arrows (i.e., logical operations) in a logical category.\par
\textbf{Conjunction.} Since conjunction is the only logical operation that is almost unchanged when transferred from traditional logic, we can define it as the following characteristic pair of arrows:\begin{center}
\includegraphics{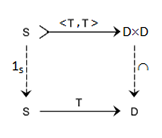}\end{center}
— i.e., as the characteristic of the product <true, true>. In this case, it is obvious that the conjunction is the greatest lower bound of its two arguments (in the order-theoretic sense). And, therefore, in the algebra of subobjects of logical categories, there exist greatest lower bounds of finite collections of subobjects (according to the lemma of pullbacks). We see that the categorical definition of conjunction is carried over almost unchanged to our generalization of logic. As a consequence of this, we also obtain the usual definition of the order of subobjects and the logical operation of implication:\par
\textbf{Implication.} The subobject of ordered pairs is an equalizer for the arrows of conjunction and the projection of the product onto the first factor. The implication is given by the characteristic functions of this subobject:\newline\begin{center}
~~~~ ~~~~ ~~ \includegraphics{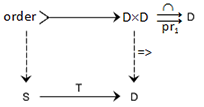}\newline\end{center}

That is — a pair of characteristic functions for the order relation induced by conjunction (when one of the subobjects is equal to the exact lower bound of its pair with another subobject).\par
\textbf{Disjunction.} And here, too, everything is analogous to the disjunction in an ordinary topos — we have the characteristic functions of embedding a direct sum into a direct product:\begin{center}
\includegraphics{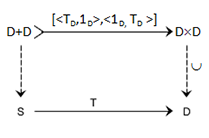}\end{center}
However (see above), we should remember that in our generalization of logic we can do without disjunction altogether (of course, at the expense of some impoverishment of the set of logical tools). The situation is more complex with negation:\par
\textbf{Negation and arrow false.} The existence of the initial object is not required here. Furthermore, there may be no minima in decreasing chains of subobjects (including the subobject \textbf{\textit{D)}}. Therefore, there is no logical value "false"; there is no negation either (in their narrow, literal sense). That is — they can only exist in the particular case of the presence of an initial object \textbf{\textit{0}}. A non-singleton generalization of logic is positive logic — logic without negation.
Thus, the general plan of proving the following theorem becomes clear:\par
\textbf{Theorem: The algebra of subobjects of any object of a logical category is a distributive lattice bounded above.}\par
Informally speaking, for a topos it is proved that this algebra is a Heyting algebra; the latter is such a lattice; moreover, in logical categories there are no negations and complements, and no specific properties of the topos are used to prove that its subobject algebras are Heyting algebras.\par
However, the impoverishment of theoretical constructions in positive logic can be avoided if we admit the existence of relative negations, when some arbitrary ("sufficiently small") subobject \textbf{\textit{D}} is declared empty (false):\begin{center}
\includegraphics{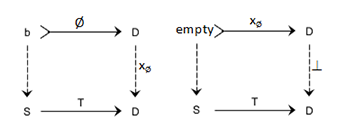}\end{center}
So, if there is a certain subobject in the classifier whose intersection with true is strictly smaller than true, then such a subobject can act as a "relatively" empty one, and its characteristic function will be a "relative" negation. And then the algebra of subobjects will become a Heyting algebra, as well as the algebra of the subobject classifier.
There may be more than one such relatively empty subobject of the classifier (depending on the category under consideration), up to an infinitely decreasing unbounded sequence (or even some collection of such sequences). In the latter case, we get a "tower" of Heyting algebras, embedded one into another as sets (these algebras will have different zeros). This will correspond to some "tower of topologies" induced by the classifier on an arbitrary algebra of subobjects of the category under consideration.
Thus, logic, impoverished by the lack of negation, will turn out to be quite enriched by the multivariance of negation.\begin{center}
\textbf{That's why:}\end{center}\par
\textbf{The algebra of subobjects of a logical category is a distributive lattice bounded above with the possibility of extending to a system of Heyting algebras (if some of the subobjects are taken as zero; Heyting algebras from this system can be included in each other as sets, but not as subalgebras — they have different zeros).}\par
From positive logic, in which there is no absolute lie, it will be natural to move to logic in which truth is also relative.\newline
\section{Projective logic}\par
But an upper bound is also not at all necessary — truth and falsity are equally ambiguous and even interchangeable: we thus come to logic as a distributive lattice of truth values. An upper bound was imposed on us by the presence of a unit arrow as the upper bound in the algebra of subobjects of a particular object. And just as we talked about the multiplicity of possible false (empty) subobjects in positive logic, we can also talk about the plurality of truths.\par
Consider the Boolean algebra B and its arbitrary interval:
$$Q_{a,b}\equiv \{x\in B\mid[a,b]\subseteq B) \wedge (a\leq x\leq b)\}$$
Let $x\in Q_{a,b}$ and $\neg x$ there is the complement of element \textbf{\textit{x}} in algebra \textbf{\textit{B}}. Then for arbitrary $y\in Q_{a,b}$\par
$(x\wedge y)\in Q_{a,b}$ and $(x\vee y)\in Q_{a,b}$; in addition, for $z\in B$, due to distributivity and $a\leq b$, we obtain:
$$(a\vee z)\wedge(a\vee b)=z=(a\wedge b)\vee(a\wedge z)$$
from where:
$$a\vee(z\wedge b)=(a\vee z)\wedge b)$$
hence:
$$((a\vee \neg x)\wedge b)\in Q_{a,b}$$
let:
$$x^{a}_{b}\equiv (a\vee\neg x\wedge b)$$
(let's call this the local complement of \textbf{\textit{x}} in the interval \textbf{[\textit{a,b}]}) this element will be the complement of \textbf{\textit{x}} in the interval $Q_{a,b}$.
$$(a\vee\neg x\wedge b)\vee x=b$$
$$(a\vee\neg x\wedge b)\wedge x=a$$
But this means that $x^{a}_{b}$ is the complement of the element \textbf{\textit{x}} in the algebra $Q_{a,b}$, which will thus be precisely a Boolean algebra (where b is the largest element and a the smallest).\par
Consequently, every interval of a Boolean algebra is itself a Boolean algebra (without being a subalgebra of \textbf{\textit{B}},— in the general case, they have different zeros, ones, and negations).\par
Moreover, the same is true for any Heyting algebra. Moreover, this is true for any modular lattice that is bounded above and below and has complements. More specifically, in the case of a Heyting algebra we have the intersection of the least upper bound of the set of elements from \textbf{\textit{B}} and \textbf{\textit{b}}, which will be the corresponding set from the interval $Q_{a,b}$; and this is Heyting negation.\par
In the case of a modular lattice with complements (no matter Boolean or Heyting), we have the same thing; in boolean form:
$$x\wedge x^{a}_{b}=a$$
$$x\vee x^{a}_{b}=b$$\par
In the Heyting case, this (as mentioned above) complement within the interval (in the general case) will be strictly less than \textbf{\textit{b}}.\par
From these considerationsit obviously follows that if some interval of a distributive (or modular) lattice has local complements, then all its subintervals also have local complements.\par
The axiom of local complementation can be expressed as follows:
$$\forall(a,b,x)|(b\geq x\geq a)\exists y(y\in Q_{a,b}):((x\wedge y)=a)  \wedge  ((x\vee y)=b)$$\par
A more rigorous presentation. In the definitions below, the symbols $\vee$ and $\wedge$ stand for upper and lower lattice faces (but when we give them a logical meaning, they, of course, become operations: “or”, “and”).\par
\textbf{\textit{Definition 2.}} For a modular (distributive) lattice \textbf{\textit{R}} and for a pair $a,b\in R$, we call the set $\{x\in R|(b\geq x\geq a)\}\equiv Q_{a,b}$ an interval; Let's call $Q_{c,d}$ a subinterval $Q_{a,b}$ if $b\geq d\geq c\geq a$;
Then $Q_{a,b}$ will be a modular (distributive) sublattice of R bounded above and below.
\textbf{\textit{Definition 3.}} An element $x^{a}_{b}\in Q_{a,b}$ is called a local \textbf{\textit{(a,b)}}-complement to $x\in Q_{a,b}$ if $x^{a}_{b}\equiv\vee\{y\in Q_{a,b}|(y\wedge x)=a\}$
Don't forget that this upper bound must belong to the interval $Q_{a,b}$; in the “Boolean” case (the quotes here mean possible modularity, but not distributivity of the lattice \textbf{\textit{R}}) we have $(x\wedge x^{a}_{b})=a$ and $(x\vee x^{a}_{b})=b$\par
\textbf{\textit{Theorem: If an interval of a modular lattice has local complements, then all its subintervals also have them. Moreover, for $(c\geq b\geq x\geq a\geq d)$ we have: $x^{a}_{b}=(a\vee x^{d}_{c}\wedge b)$.}}\par
Proof. Let $Q^{a}_{b}$ be a subinterval of an interval with local complements and let $x\in Q_{a,b}$ have such local complement $\neg x$; then, due to modularity and $b\geq x$, we obtain: $a\vee(x\wedge b)=(a\vee x)\wedge b$ and:
$$((a\vee\neg x)\wedge b)\wedge x=(a\vee\neg x)\wedge x=a$$\par
Moreover, $x\vee (a\vee (\neg x)\wedge b)=x\vee(\neg x\wedge b)$, as the supremum of the wider set, after intersecting with \textbf{\textit{b}}, will be the corresponding maximal element.\par
The same is true for the stronger case of a distributive lattice. We have the double projection operator into the interval as a lattice homomorphism: $x\longrightarrow(a\vee x\wedge b)$. In the case of Boolean algebra, all the proofs are generally quite simple — it is enough to take the corresponding upper and lower bounds. The proof is over.\newline\par
\textbf{\textit{Corollary: As a result of projecting an interval of a Heyting (or Boolean) algebra, the interval also inherits the structure of a Heyting (or Boolean) algebra.}}\newline\par
These structures, however, will not be subalgebras—complements, maximal and minimal elements are not preserved.
General conclusion. The most general model of logic is a modular lattice that is unbounded above and below and has local complements. We call this logic projective (by analogy with projective geometry).\par
\textbf{\textit{Comment.}} In a Heyting algebra, the least upper bound of an arbitrary (not necessarily finite) set of elements is taken; since we are talking about constructive negation, we get violations of the Boolean law of excluded middle; but quantum logic also tells us about violations of the Boolean law of non-contradiction; Therefore, we need to allow the possibility of considering any precise bounds of any subsets of the lattice — which would be a further generalization of logical algebras — with different models of implication and negation, as well as the return of Boolean symmetry of true and false. Those. — projective logic should include the possibility of inverting the lattice order. But the main idea here is the relativity of truth, falsity and negation — this is what we wanted to substantiate in the article.\newline
\section{A note on the semantics of Galois correspondences}\par
In conclusion, a few words about logical semantics. We take an individual region open from above—i.e. a set of elements from a certain universal set, which we do not consider (we deliberately do not talk about sets, but only about aggregates, because we try to apply only some of the axioms of set theory). We declare some of these elements to be predicates. We believe that sets of elements (some of which are predicates) can form Cartesian powers. We establish the following correspondence.\par
To each of the predicates we associate the set of all sequences of elements on which it is true (in this case, predicates can also be elements of sequences) and to each set of n-sequences of elements we associate the set of all predicates that are true on it:
$$P^{*}\equiv \{x_{i(i\in 1...n)}\mid P(x_{i})\}$$
$$X^{*}\equiv \{P\mid P(x_{i})\}$$\par
This will be a Galois correspondence (by definition). In such a construction one can introduce arbitrary operations on elements or predicates that preserve the truth relation. Then the closure systems induced by this Galois correspondence on predicates and on sequences of elements:
$$P\longrightarrow P^{**}$$
$$X\longrightarrow X^{**}$$
— will set the derivable formulas for a theory that is true on a certain set of sequences of elements. Such closure systems can be algebraic or topological. It is fairly obvious to interpret quantifiers in terms of this Galois correspondence. We write the axioms of the theory. This is how we get a logic in which individuals and properties are not only relative, but also on an equal footing. Of course, this rough picture needs further detail.\newline\newline\newline
\section{References}
\noindent [1]  Haskell B. Curry, McGraw-Hill Book Company, INC.\newline
\textit{Foundations of Manthematical logic}, (1969).\newline
[2] Garrett Birkhoff, Providence Rhode Island,\newline
\textit{Lattice theory}, (1967).\newline
[3] Helena Rasiowa and Roman Sikorski,  Panstwowe Wydawnlctwo Naukowe Warszawa,\newline
 \textit{The Mathemayics of Metamathematics}, (1963).\newline
 [4] R. Goldblatt. North-Holland Publishing Company: Amsterdam, New York, Oxford.\newline
 \textit{Topoi. The categorial analysis of logic}, (1979).\newline
 [5] P.M. Kon,  D. Reidel Publishing Company,\newline
 \textit{Universal Algebra}, (1961).\newline
 [6] A. Heyting, North-Holland publishing company, Amsterdam,\newline
 \textit{Intuitionism, An introduction}, (1956)
\end{document}